\documentclass{amsart}
\usepackage{amssymb} 
\usepackage{color} 

\newtheorem{theorem}{Theorem}
\newtheorem{lemma}{Lemma}

\newtheorem{remark}{Remark}
\newcommand{\pp}{\noindent {\em Proof. }}
\begin{document}
\title{Identities of graded simple algebras}

\author[D. D. Repov\v s and M. V. Zaicev]
{Du\v san Repov\v s and Mikhail Zaicev}

\address{Du\v san D. Repov\v s \\Faculty of Education, and
Faculty of  Mathematics and Physics, University of Ljubljana,
 Ljubljana, 1000, Slovenia}
\email{dusan.repovs@guest.arnes.si}

\address{Mikhail V. Zaicev \\Department of Algebra\\ Faculty of Mathematics and
Mechanics\\  Moscow State University \\ Moscow,119992, Russia}

\email{zaicevmv@mail.ru}

\thanks{The first author was supported by the SRA
grants P1-0292-0101, J1-7025-0101, J1-6721-0101 and J1-5435-0101. 
The second author was partially supported 
by the RFBR grant 16-01-00113a. We thank the referee for comments and suggestions.}

\keywords{Polynomial identities, graded algebras, codimensions,
exponential growth}

\subjclass[2010]{Primary 17B01, 16P90; Secondary 16R10}

\begin{abstract}
We study identities of finite dimensional algebras over a field of 
characteristic zero, graded by an arbitrary groupoid $\Gamma$. First we prove 
that its graded colength has a polynomially bounded growth. For any graded
simple algebra $A$ we prove the existence of the graded PI-exponent, provided that
$\Gamma$ is a commutative semigroup. If $A$ is simple in a non-graded sense
the existence of the graded PI-exponent is proved without any restrictions on $\Gamma$.
\end{abstract}

\date{\today}

\maketitle
\vskip 0.2in

\section{Introduction}

We study numerical characteristics of identities of finite dimensional graded
simple algebras over a field of characteristics zero. The main object of our investigations 
is the asymptotic behaviour of sequences of graded codimensions and graded colengths of 
such algebras (all necessary definitions and notions will be given in the next 
section). Given a graded algebra $A$, one can associate the sequence of so-called {\it graded
codimensions} $\{c_n^{gr}(A), n=1,2,\ldots\}$. This sequence is an important numerical 
invariant of graded identities of $A$. It is known that this sequence is exponentially
bounded, that is  $c_n^{gr}(A) \le a^n$ for some real $a$, provided that $\dim A<\infty$.
In this case the following natural question arises: does the limit
\begin{equation}\label{int1}
exp^{gr}(A)=\lim_{n\to\infty} \sqrt[n]{c_n^{gr}(A)}
\end{equation}
exist and what are its possible values? If the limit (\ref{int1}) exists then it is called
the {\it graded PI-exponent} of $A$. 

In the non-graded case, codimension growth is well understood. Existence and integrality of
the (non-graded) PI-exponent was conjectured by Amitsur in 1980's for associative PI-algebras.
Amitsur's conjecture was confirmed in \cite{GZAdv98, GZAdv99}. Later the same result 
was proved for finite dimensional Lie algebras \cite{GRZ1,GRZ2,Z-Izv},
Jordan and alternative algebras \cite{GZLaa2007,GZTams2010,GShZ2011} and
many other algebraic systems. In the general nonassociative case, for any real $\alpha>1$,
examples of algebras with PI-exponent equal to $\alpha$ were constructed in \cite{GMZTams}.
Recently, the first example of algebra $A$ such that the PI-exponent of $A$ does not
exist, was constructed \cite{Z-ERA}. Nevertheless, for any finite dimensional simple algebra the PI-exponent
does exist \cite{GZ-JLMS}.

For graded algebras there are only partial results of this kind. For example, if $A$ is
an associative graded PI-algebra then its graded PI-exponent always exists and it is an
integer \cite{GLm}. An existence of the ${\mathbb Z}_2$-graded PI-exponent for any finite 
dimensional simple Lie superalgebra has recently been proved  in \cite{RZ-Ark}. Note that for
finite dimensional Lie superalgebras, both graded and ordinary PI-exponents can be
fractional \cite{GZ-JLMS, RZ, RZ-b2}. The main purpose of this paper is to prove 
the existence of graded PI-exponents for any finite dimensional graded simple algebra
(see Theorem \ref{t2}).

Another important numerical characteristic of identities of an algebra $A$ is the  so-called 
{\it colength sequence} $l_n(A)$. Except for its independent interest, asymptotic behaviour of
$\{l_n(A)\}$ plays an important role in the studies of asymptotics of $\{c_n(A)\}$. The polynomial 
type of growth of $\{l_n(A)\}$ is very convenient for investigations of codimension
growth.

Polynomial upper bounds of the colength for any associative PI-algebra were established in \cite{Be}.
For an arbitrary (non-associative) finite dimensional algebra the same restriction was
obtained in \cite{GMZ-AAM}. In the case of finite dimensional Lie superalgebras, polynomial growth 
of ${\mathbb Z}_2$-graded colength has recently been confirmed  in \cite{Z-VMSU}. In order to get 
the main result of the paper we will find the polynomial upper bound for graded colength of a
finite dimensional graded algebra (see Theorem \ref{t1}).
\vskip 0.2in

\section{Preliminaries}

Let $\Gamma$ be a groupoid. An $F$-algebra $A$ is said to be $\Gamma$-{\it graded} if there is a vector 
space decomposition
$$
A=\bigoplus_{g\in\Gamma} A_g
$$
and $A_g A_h\subseteq A_{gh}$ for all $g,h\in \Gamma$. An element $a\in A$ is called {\it homogeneous
of degree} $g$ if $a\in A_g$ and  in this case we write $\deg_\Gamma a=g$. A subspace  $V\subseteq A$
is homogeneous iff $V=\bigoplus_{g\in \Gamma}(V\cap A_g)$. We call $A$ {\it graded simple} if it has no
homogeneous ideals. For instance, if $\Gamma$ is a group and $A=F[\Gamma]$ is its group algebra then
$A$ is $\Gamma$-graded simple but is not simple in the usual sense. On the other hand, any simple 
algebra with an arbitrary grading is graded simple.

We recall some key notions from the theory of graded and ordinary identities and their numerical invariants. 
We refer the reader to \cite{GZbook, DRbook} for details. Consider an absolutely free algebra
$F\{X\}$ with a free generating set
$$
X=\bigcup_{g\in\Gamma} X_g,\quad |X_g|=\infty \quad \mbox{for any} \quad g\in\Gamma.
$$
One can define a $\Gamma$-grading on $F\{X\}$ by setting $\deg_\Gamma x=g$, when $x\in X_g$, and
extend this grading to the entire $F\{X\}$ in the natural way. A polynomial $f(x_1,\ldots, x_n)$
in homogeneous variables $x_1\in X_{g_1},\ldots, x_n\in X_{g_n}$ is called a {\it graded identity}
of a $\Gamma$-graded algebra $A$ if $f(a_1,\ldots, a_n)=0$ for any  
$a_1\in A_{g_1},\ldots, a_n\in A_{g_n}$. The set $Id^{gr}(A)$ of all graded identities of $A$
forms an ideal of $F\{X\}$ which is stable under graded homomorphisms $F\{X\}\rightarrow F\{X\}$.

First, let $\Gamma$ be finite, $\Gamma=\{g_1,\ldots,g_t\}$ and 
$X=X_{g_1}\bigcup\ldots\bigcup X_{g_t}$. Denote by $P_{n_1,\ldots,n_t}$ the subspace of  $F\{X\}$
of multilinear polynomials of total degree $n=n_1+\cdots+n_t$ in variables
$x^{(1)}_{1},\ldots,x^{(1)}_{n_1}\in X_{g_1},\ldots, x^{(t)}_{1},\ldots,x^{(t)}_{n_t}\in X_{g_t}$.
Then the value
$$
c_{n_1,\ldots,n_t}(A)=\dim\frac{P_{n_1,\ldots,n_t}}{P_{n_1,\ldots,n_t}\cap Id^{gr}(A)}
$$
is called a {\it partial codimension} of $A$ while
\begin{equation}\label{e1}
c_n^{gr}(A)=\sum_{n_1+\cdots+n_t=n} {n\choose n_1,\ldots,n_t} c_{n_1,\ldots,n_t}(A)
\end{equation}
is called a {\it graded codimension} of $A$. Recall that the {\it support of the grading} is the set
$$
Supp~A=\{g\in \Gamma\vert A_g\ne 0\}.
$$
Note that if $Supp~A\ne\Gamma$, say,  $Supp~A=\{g_1,\ldots,g_k\}$, $k<t$, then the value
\begin{equation}\label{e2}
\sum_{n_1+\cdots+n_k=n} {n\choose n_1,\ldots,n_k} \dim\frac{P_{n_1,\ldots,n_k}}{P_{n_1,\ldots,n_k}\cap Id^{gr}(A)}
\end{equation}
coincides with (\ref{e1}). This allows us to consider (\ref{e2}) as the definition of the graded codimension
of $A$ even if $\Gamma$ is infinite, provided that $Supp~A=\{g_1,\ldots,g_k\}$. 

For convenience, denote
\begin{equation}\label{e3}
P_{n_1,\ldots,n_k}(A)=\frac{P_{n_1,\ldots,n_k}}{P_{n_1,\ldots,n_k}\cap Id^{gr}(A)}.
\end{equation}
Given $1\le j\le k$, consider the action of the symmetric group $S_{n_j}$ on $P_{n_1,\ldots,n_k}$
defined by
$$
\sigma f(\ldots,x_1^{(j)},\ldots, x_{n_j}^{(j)},\ldots) = 
f(\ldots,x_{\sigma(1)}^{(j)},\ldots, x_{\sigma(n_j)}^{(j)},\ldots).
$$
Then the spaces $P_{n_1,\ldots,n_k}$ and $P_{n_1,\ldots,n_k}(A)$ become $F[H]$-modules, where
$H=S_{n_1}\times\cdots\times S_{n_k}$. Any $F[H]$-module $P_{n_1,\ldots,n_k}(A)$ is decomposed
into the sum of irreducible $F[H]$-submodules and in the languages of group characters it can be
written as
\begin{equation}\label{e4}
\chi_{_H}(P_{n_1,\ldots,n_k}(A))=\sum_{\lambda^{(1)}\vdash n_{1},\ldots,\lambda^{(k)}\vdash n_{k}}
m_{\lambda^{(1)},\ldots,\lambda^{(k)}} \chi_{_{\lambda^{(1)},\ldots,\lambda^{(k)}}}.
\end{equation}
Here, $\chi_{_{\lambda^{(1)},\ldots,\lambda^{(k)}}}$ is the character of the irreducible $H$-representation 
defined by the $k$-tuple $(\lambda^{(1)},\ldots,\lambda^{(k)})$ of partitions
$\lambda^{(1)}\vdash n_{1},\ldots,\lambda^{(k)}\vdash n_{k}$ and
$m_{\lambda^{(1)},\ldots,\lambda^{(k)}}$ is the multiplicity of the corresponding $F[H]$-module
in $P_{n_1,\ldots,n_k}(A)$. The integer
\begin{equation}\label{e5}
l_{\lambda^{(1)},\ldots,\lambda^{(k)}}(A) =
\sum_{\lambda^{(1)}\vdash n_{1},\ldots,\lambda^{(k)}\vdash n_{k}}  m_{\lambda^{(1)},\ldots,\lambda^{(k)}}
\end{equation}
is called the {\it partial colength}, whereas the integer
\begin{equation}\label{e6}
l_{n}^{gr}(A) =
\sum_{n_1+\cdots+n_k=n}  l_{n_1,\ldots,n_k}(A)
\end{equation}
is called the {\it graded colength} of $A$.

As it was mentioned in the introduction, graded codimensions are exponentially bounded if 
$A$ is finite dimensional. Namely,
\begin{equation}\label{e7}
c_{n}^{gr}(A) \le d^{n+1}
\end{equation}
where $d=\dim A$ (see \cite{BD-Laa} and also \cite[Proposition 2]{GZTams2010}). This
result was proved in \cite{GZTams2010, BD-Laa}  under the assumption that $\Gamma$
is a finite group. The same argument is valid for an arbitrary groupoid. Relation
(\ref{e7}) allows us to consider upper and lower limits of $\sqrt[n]{c_n^{gr}(A)}$ and
we can define the lower and the upper graded PI-exponents as follows:
$$
\underline{exp}^{gr}(A)=\liminf_{n\to\infty} \sqrt[n]{c_n^{gr}(A)},\qquad
\overline{exp}^{gr}(A)=\limsup_{n\to\infty} \sqrt[n]{c_n^{gr}(A)}.
$$
If the lower and the upper limits coincide then we also define the graded PI-exponent by
$$
exp^{gr}(A)=\underline{exp}^{gr}(A)=\overline{exp}^{gr}(A).
$$

Representation theory of symmetric groups is a useful tool for studying asymptotics
of codimension growth. Basic notions of $S_n$-representations can be found  in
\cite{JK} and its application to PI-theory  in \cite{GZbook, DRbook}.

Recall that, given a partition $\lambda\vdash n$, there is exactly one (up to
isomorphism)  irreducible $S_n$-representation defined by $\lambda$. Its
character and dimension are denoted by $\chi_{_\lambda}$ and $\chi_{_\lambda}(1)=d_\lambda$,
respectively. For the group $H=S_{n_1}\times\cdots\times S_{n_k}$ any
irreducible representation is defined by the $k$-tuple of partitions
$\lambda^{(1)}\vdash n_{1},\ldots,\lambda^{(k)}\vdash n_{k}$ and its character 
and dimension are $\chi_{_{\lambda^{(1)},\ldots,\lambda^{(k)}}}$. Moreover,
\begin{equation}\label{e8}
\chi_{_{\lambda^{(1)},\ldots,\lambda^{(k)}}}(1)=d_{\lambda^{(1)}}\cdots 
d_{\lambda^{(k)}},
\end{equation}
respectively. In particular, (\ref{e4}) and (\ref{e8}) imply the equality
\begin{equation}\label{e9}
c_{n_1,\ldots,n_k}(A)=\chi_{_H}(P_{n_1,\ldots,n_k}(A))(1)= 
\end{equation}
$$
\sum_{\lambda^{(1)}\vdash n_{1},\ldots,\lambda^{(k)}\vdash n_{k}}
m_{\lambda^{(1)},\ldots,\lambda^{(k)}} d_{\lambda^{(1)}}\cdots d_{\lambda^{(k)}}.
$$

Let $d\ge 1$ be a fixed integer and let $\nu=(\nu_1,\ldots,\nu_q)\vdash m$ be a
partition of $m$ with $q\le d$. Dimension of an irreducible $F[S_m]$-module
with the character $\chi_{_\nu}$ is closely connected with the following function:
$$
\Phi(\nu) = \frac{1}{(\frac{\nu_1}{m})^\frac{\nu_1}{m}\cdots (\frac{\nu_d}{m})^\frac{\nu_d}{m}}.
$$
Here we assume that $\nu_{q+1}=\ldots=\nu_d=0$ in the case $q<d$ and $0^0=1$. The
values $\Phi(\nu)^m$ and $d_\nu$ are close in the following sense.
\vskip 0.2in

\begin{lemma}\label{la1}({\rm see} \cite[Lemma 1]{GZ-JLMS})
Let $m\ge 100$. Then
$$
\frac{\Phi(\nu)^m}{m^{d^2+d}} \le d_\nu \le m \Phi(\nu)^m.
$$
\end{lemma}
\hfill $\Box$

We will use the following property of $\Phi$. Let $\nu$ and $\rho$ be any two partitions 
of $m$, such that $\nu=(\nu_1,\ldots,\nu_p)$, $\rho=(\rho_1,\ldots,\rho_q)$, $p,q\le d$ and
$q=p$ or $q=p+1, \rho_{p+1}=1$. As before, we consider $\rho$ and $\nu$ as partitions
with $d$ components. We say that the Young diagram $D_\rho$ is obtained from diagram $D_\nu$
by pushing down one box if there exist $1\le i<j\le d$ such that $\rho_i=\nu_i-1$, 
$\rho_j=\nu_j+1$ and $\rho_t=\nu_t$ for all remaining $1\le t\le d$.
\vskip 0.2in

\begin{lemma}\label{la2} ({\rm see} \cite[Lemma 3]{GZ-JLMS}, \cite[Lemma 2]{ZR})
Let $D_\rho$ be obtained from $D_\nu$ by pushing down one box. Then 
$\Phi(\rho) \ge \Phi(\nu)$.
\end{lemma}
\hfill $\Box$
\vskip 0.2in

\section{Polynomial growth of graded colength}

Consider a finite dimensional $\Gamma$-graded algebra $A$ with the support
$Supp~A=\{g_1,\ldots,g_k\}$, $A=A_{g_1}\oplus\cdots\oplus A_{g_k}$. Let
$$
d_1=\dim A_{g_1},\ldots, d_k=\dim A_{g_k}
$$
be dimensions of the homogeneous components. Recall that an irreducible $F[S_t]$-module
corresponding to the partition $\mu\vdash t$ can be realized as a minimal left
$F[S_t]$-ideal generated by an essential idempotent $e_{T_\lambda}$ where $T_{\lambda}$
is some Young tableaux with Young diagram $D_\lambda$. For 
$H=S_{n_1}\times\cdots\times S_{n_k}$, any irreducible $F[H]$-module is isomorphic to the
tensor product of $F[S_{n_1}],\ldots,F[S_{n_k}]$-modules with characters 
$\chi_{_{\lambda^{(1)}}},\ldots,\chi_{_{\lambda^{(k)}}}$, respectively. The following remark 
easily follows from the construction of essential idempotents and therefore we omit the proof.
\vskip 0.2in

\begin{lemma}\label{l1}
Let $\lambda^{(1)}=(\lambda^{(1)}_1,\ldots, \lambda^{(1)}_{q_1}),\ldots,
\lambda^{(k)}=(\lambda^{(k)}_1,\ldots, \lambda^{(k)}_{q_k})$ be partitions of
$n_1,\ldots, n_k$, respectively. Suppose that the multiplicity $m_{\lambda^{(1)}, 
\ldots, \lambda^{(k)}}$ on the right hand side of (\ref{e4}) is nonzero. Then
$q_1\le d_1,\ldots, q_k\le d_k$.
\end{lemma}
\hfill $\Box$

For convenience we shall assume as before that $q_1= d_1,\ldots, q_k= d_k$ even if $q_i$ is 
strictly less than $d_i$ for some $i$.

Denote by $R=R\{X_{g_1}\cup\ldots\cup X_{g_k}\}$ the relatively free algebra of the variety
$var~A$ of $\Gamma$-graded algebras generated by $A$. Denote by
$R^{n_1,\ldots, n_k}_{d_1,\ldots, d_k}$ the subspace of polynomials in $R$ of degree
$n_1$ in the set of variables $\{X_1^{(1)},\ldots,X_{d_1}^{(1)}\}\subseteq X_{g_1}$, of
degree $n_2$ in $\{X_1^{(2)},\ldots,X_{d_2}^{(2)}\}\subseteq X_{g_2}$, etc.
\vskip 0.2in

\begin{lemma}\label{l2}
Multiplicities on the right hand side of (\ref{e4}) satisfy the inequalities
$$
m_{\lambda^{(1)},\ldots,\lambda^{(k)}} \le \dim R^{n_1,\ldots, n_k}_{d_1,\ldots, d_k}.
$$
\end{lemma}
\pp
Let $\widetilde P_{n_1,\ldots,n_k}$ be the subspace of multilinear polynomials of degree 
$n=n_1+\cdots+n_k$ on $x^{(1)}_1,\ldots, x^{(1)}_{n_1},\ldots, 
x^{(k)}_1,\ldots, x^{(k)}_{n_k}$ in $R$. Here $x^{(i)}_j\in X_{g_i}$ for all $1\le i\le k,
1\le j\le n_i$. Then $\widetilde P_{n_1,\ldots,n_k}$ is isomorphic to $P_{n_1,\ldots,n_k}(A)$
as an $F[H]$-module. Denote for brevity $q=m_{\lambda^{(1)},\ldots,\lambda^{(k)}}$ and
consider the $F[H]$-submodule
$$
M=M_1\oplus\cdots\oplus M_q
$$
of $\widetilde P_{n_1,\ldots,n_k}$, where $M_1,\ldots, M_q$ are isomorphic irreducible
$F[H]$-modules with $H$-character $\chi_{_{\lambda^{(1)},\ldots,\lambda^{(k)}}}$. Any
$M_j$ is generated as an $F[H]$-module by a multilinear polynomial of the type
$$
f_j(x^{(1)}_1,\ldots, x^{(1)}_{n_1},\ldots,x^{(k)}_1,\ldots, x^{(k)}_{n_k})=
$$
$$
 e_{T_{\lambda^{(1)}}}\cdots e_{T_{\lambda^{(k)}}}
 h_j(x^{(1)}_1,\ldots, x^{(1)}_{n_1},\ldots,x^{(k)}_1,\ldots, x^{(k)}_{n_k})
$$
with a multilinear polynomial $h_j\in \widetilde P_{n_1,\ldots,n_k}$.

One can split the set of indeterminates $x^{(1)}_1,\ldots, x^{(k)}_{n_k}$ into a 
disjoint union of subsets
$$
P^{(1)}_1= \{x^{(1)}_1,\ldots, x^{(1)}_{\lambda^{(1)}_1}\}, 
$$
$$
P^{(1)}_2= \{x^{(1)}_{\lambda^{(1)}_1+1},\ldots, x^{(1)}_{\lambda^{(1)}_1+\lambda^{(1)}_2}\},
$$
$$
\cdots
$$
$$
P^{(1)}_{d_1}= \{x^{(1)}_{n_1-\lambda^{(1)}_{d_1}+1},\ldots, x^{(1)}_{n_1}\},
$$
$$
\ldots
$$
$$
P^{(k)}_1= \{x^{(k)}_1,\ldots, x^{(k)}_{\lambda^{(k)}_1}\}, 
$$
$$
P^{(k)}_2= \{x^{(k)}_{\lambda^{(k)}_1+1},\ldots, x^{(k)}_{\lambda^{(k)}_1+\lambda^{(k)}_2}\},
$$
$$
\cdots
$$
$$
P^{(k)}_{d_k}= \{x^{(k)}_{n_k-\lambda^{(1)}_{d_k}+1},\ldots, x^{(k)}_{n_k}\}
$$
so that all $f_1,\ldots,f_q$  are symmetric on any subset $P^{(1)}_{1},\ldots,
P^{(k)}_{d_k}$.

Now we identify all variables in each symmetric subset, that is, we apply a homomorphism
$\varphi: R\to R$ such that
\begin{itemize}
\item
$\varphi(x^{(j)}_\alpha)= x^{(j)}_{1}$ if $1\le\alpha\le\lambda^{(j)}_1$,
\item[]$\cdots$
\item
$\varphi(x^{(j)}_\alpha)= x^{(j)}_{d_j}$ if
$\lambda^{(j)}_1+\cdots+\lambda^{(j)}_{d_j-1}< \alpha\le n_j$
\end{itemize}
for all $j=1,\ldots,k$. Then all $\varphi(f_1),\ldots,\varphi(f_q)$ lie in
$R^{n_1,\ldots,n_k}_{d_1,\ldots,d_k}$. Note that the total linearization of
each of $\varphi(f_j)$ is equal to $f_j$ with a nonzero coefficient independent of $j$.
Hence a nontrivial linear relation $\alpha_1\varphi(f_1)+\cdots+\alpha_q\varphi(f_q)=0$
implies the same relation $\alpha_1 f_1+\cdots+\alpha_q f_q=0$. But $f_1,\ldots,f_q$ 
belong to distinct irreducible summands $M_1,\ldots, M_q$, respectively. In particular,
they are linearly independent. Hence $q$ does not exceed $\dim R^{n_1,\ldots,n_k}_{d_1,\ldots,d_k}$,
and the proof is completed.
\hfill $\Box$

Now we restrict the dimension of $R^{n_1,\ldots,n_k}_{d_1,\ldots,d_k}$.
\vskip 0.2in

\begin{lemma}\label{l3}
Let $A=A_{g_1}\oplus\cdots\oplus A_{g_k}$ be a $\Gamma$-graded algebra with the support
$\{g_1,\ldots,g_k\}$ and let $d_1=\dim A_{g_1},\ldots, d_k=\dim A_{g_k}$. Then
\begin{equation}\label{e10}
\dim R^{n_1,\ldots,n_k}_{d_1,\ldots,d_k} \le (d_1+\cdots+d_k)(n_1+1)^{d_1^2}\cdots (n_k+1)^{d_k^2}.
\end{equation}
\end{lemma}
\pp
Let $\{a^{(g_i)}_1,\ldots, a^{(g_i)}_{d_i}\}$ be a basis of the subspace $A_{g_i}, 1\le i\le k$.
Consider a polynomial ring $F[Y]$, where $Y=Y_{g_1}\cup\ldots\cup Y_{g_k}$ and
$$
Y_{g_i}=\{ y^{g_i}_{m,j}|~1\le m\le d_i,j=1,2,\ldots \}.
$$
Then algebra $F[Y]$ can be naturally endowed by a $\Gamma$-grading with $Supp~F[Y]=\{g_1,\ldots,g_k\}$
if we set $\deg_\Gamma y=g_i$ when $y\in Y_{g_i}$. Denote $\widetilde A = A\otimes F[Y]$
and fix elements
$$
z^{g_i}_j = \sum_{m=1}^{d_i} a^{(g_i)}_m\otimes y^{g_i}_{m,j},\quad j=1,2,\ldots,
$$
in $\widetilde A$. Then $alg\{z^{g_i}_j\}$ is also a $\Gamma$-graded algebra, where $\deg z^{g_i}_j = g_i$.
Moreover, $\widetilde A\simeq R$ and $R^{n_1,\ldots,n_k}_{d_1,\ldots,d_k}$ is a subspace of
$A\otimes T$, where $T$ is the subspace of $F[Y]$ spanned by monomials of degree at most $n_t$ in the
set of indeterminates $\{ y^{g_t}_{m,j}|~1\le m,j\le d_t \}$, $t=1,\ldots,k$. Clearly,
\begin{equation}\label{e11}
\dim T \le (n_1+1)^{d_1^2}\cdots (n_k+1)^{d_k^2}
\end{equation}
hence (\ref{e10}) follows from (\ref{e11}).
\hfill $\Box$

Now we are ready to get an upper bound for graded colength of $A$.
\vskip 0.2in

\begin{theorem}\label{t1}
Let $A=\bigoplus_{g\in\Gamma} A_g$ be a finite dimensional algebra graded by groupoid $\Gamma$
with $Supp~A=\{g_1,\ldots, g_k\}$. Let also $\dim A_{g_i}=d_i, 1\le  i \le k$. Then the $n$th
graded colength of $A$ satisfies the inequality
$$
l_n^{gr} \le d(n+1)^{k+d_1^2\cdots+d_k^2+d_1+\cdots+d_k}
$$
where $d=\dim A=d_1+\cdots+d_k$.
\end{theorem}

\pp
By Lemma \ref{l1}, the total number of partitions $\lambda_i\vdash n_i$ does not exceed $(n_i+1)^{d_i}$
for any $i=1,\ldots,k$. Hence, by (\ref{e5}) and Lemmas \ref{l2} and \ref{l3}, we have
$$
l_{n_1,\ldots,n_k}(A) \le d(n_1+1)^{d_1^2+d_1}\cdots (n_k+1)^{d_k^2+d_k}
$$
and
$$
l_n^{gr} \le d(n+1)^{k+d_1^2\cdots+d_k^2+d_1+\cdots+d_k}
$$
as follows from (\ref{e6}).
\hfill $\Box$
\vskip 0.2in

\section{Existence of graded PI-exponents}

We begin this section with a technical result connecting dimensions of 
irreducible representations of symmetric groups and multinomial coefficients.
Given a partition $\mu=(\mu_1,\ldots,\mu_t)\vdash m$, we denote by 
$q\mu(q\mu_1,\ldots,q\mu_t)$ the partition of $qm$, where $q\ge 1~$ is an arbitrary
integer. We also define the {\it height} $ht(\mu)$ as $ht(\mu)=t$. Recall that
$d_\mu=\chi_{_\mu(1)}$ is the {\it dimension} of the corresponding irreducible representation
of $S_m$.
\vskip 0.2in

\begin{lemma}\label{le1}
Let $n_1,\ldots,n_k$ be positive integers, $n_1+\cdots+n_k=n\ge 100$. Let also
$\lambda^{(1)},\ldots,\lambda^{(k)}$ be partitions of $n_1,\ldots,n_k$, respectively,
such that $ht(\lambda^{(1)}),\ldots,ht(\lambda^{(k)}) \le d$. If $q\ge 100$ then
$$
{qn\choose qn_1,\ldots,qn_k}d_{q\lambda^{(1)}}\cdots d_{q\lambda^{(k)}} \ge
\left(\frac{1}{qn} \right)^{k(d^2+d+1)} \left[\frac{1}{n^{2k}}
{n\choose n_1,\ldots,n_k} d_{\lambda^{(1)}}\cdots d_{\lambda^{(k)}}
\right]^q.
$$
\end{lemma}

\pp
Given nonnegative real $\alpha_1,\ldots\alpha_k$ with $\alpha_1+\cdots+\alpha_k=1$,
we denote
$$
\Phi(\alpha_1,\ldots\alpha_k)=\frac{1}{(\alpha_1)^{\alpha_1}\cdots (\alpha_k)^{\alpha_k}}~.
$$
From the Stirling formula for factorials it easily follows that
\begin{equation}\label{eq1}
\frac{1}{m^k}\Phi\left(\frac{m_1}{m},\cdots,\frac{m_k}{m}\right)^m \le
{m\choose m_1,\ldots,m_k} \le m\Phi\left(\frac{m_1}{m},\cdots,\frac{m_k}{m}\right)^m,
\end{equation}
where $m_1,\ldots,m_k$ are nonnegative integers and $m_1+\cdots+m_k=m$. Applying (\ref{eq1}) to
$$
P={qn\choose qn_1,\ldots,qn_k}
$$
we obtain
$$
P> \left(\frac{1}{qn} \right)^{k} \Phi\left(\frac{qn_1}{qn},\ldots,\frac{qn_k}{qn}\right)^{qn}
= \left(\frac{1}{qn} \right)^{k} \left[\Phi\left(\frac{n_1}{n},\ldots,\frac{n_k}{n}\right)^{n}
\right]^q.
$$
Applying again (\ref{eq1}) we get
\begin{equation}\label{eq2}
P>
\left(\frac{1}{qn} \right)^{k} \left[\frac{1}{n}{n\choose n_1,\ldots,n_k}
\right]^q.
\end{equation}
It follows from Lemma \ref{la1} and (\ref{eq1}) that
$$
d_{q\lambda^{(1)}}\cdots d_{q\lambda^{(k)}}> \left(\frac{1}{qn_1}\cdots \frac{1}{qn_k}\right)^{k(d^2+d)}
\left[\Phi(\lambda^{(1)})^{n_1}\cdots \Phi(\lambda^{(k)})^{n_k}\right]^q >
$$
\begin{equation}\label{eq3}
\left(\frac{1}{qn}\right)^{k(d^2+d)}
\left[d_{\lambda^{(1)}}\cdots d_{\lambda^{(k)}} \right]^q.
\end{equation}

Now our lemma is a consequence of (\ref{eq2}) and (\ref{eq3}).
\hfill $\Box$

Recall that $A$ is a $d$-dimensional $\Gamma$-graded algebra. Now let
\begin{equation}\label{eqq}
a=\overline{exp}^{gr}(A)=\limsup_{n\to\infty} \sqrt[n]{c_n^{gr}(A)}.
\end{equation}
The next lemma is the main step of the proof of Theorem \ref{t2}.
\vskip 0.2in

\begin{lemma}\label{le2}
Let $\Gamma$ be a commutative semigroup and let $a$ in (\ref{eqq}) be strictly greater
than $1$. If $A$ is graded simple then for any $\varepsilon>0$ and any $\delta>0$ there exists an increasing
sequence of positive integers $n^{(1)}, n^{(2)},\ldots$
such that
\begin{itemize}
\item[(i)]
$\sqrt[n]{c_n^{gr}(A)}>(1+\delta)(1-\varepsilon)$ for all $n= n^{(1)}, n^{(2)},\ldots$; and
\item[(ii)]
$n^{(q+1)}, n^{(q)}\le d$, for all $q=1,2,\ldots~~$.
\end{itemize}
\end{lemma}

\pp
Clearly, there exists an integer $n^{(1)}$ such that
$$
c_{n^{(1)}}^{gr}(A)>(a-\varepsilon)^{n^{(1)}}
$$
and $n^{(1)}$ can be chosen arbitrary large. There are also $n_1,\ldots,n_k\ge 0$ such that
$n_1+\cdots+n_k=n^{(1)}$ and
$$
{n^{(1)}\choose n_1,\ldots,n_k}c_{n_1,\ldots,n_k}(A) >
\frac{1}{(n^{(1)}+1)^k} (a-\varepsilon)^{n^{(1)}}.
$$
(see (\ref{e1})). Without loss of generality, we can suppose that $k=|Supp~A|$. Consider the
$H=S_{n_1}\times\cdots\times S_{n_k}$-action on $P_{n_1,\ldots,n_k}$. It follows from (\ref{e5}),
(\ref{e6}), (\ref{e9}) that there exist partitions $\lambda^{(1)}\vdash n_{1},
\ldots, \lambda^{(k)}\vdash n_{k}$ such that
$$
d_{\lambda^{(1)}}\cdots d_{\lambda^{(k)}}>\frac{1}{l_{n^{(1)}}^{gr}(A)} c_{n_1,\ldots,n_k}(A).
$$
By Theorem \ref{t1} we have
$$
l_{n^{(1)}}^{gr}(A)< d(n^{(1)}+1)^{k(d+1)^2}
$$
hence
$$
d_{\lambda^{(1)}}\cdots d_{\lambda^{(k)}}>\frac{1}{d(n^{(1)}+1)^{k(d+1)^2}} c_{n_1,\ldots,n_k}(A).
$$
and 
\begin{equation}\label{eq4}
{n^{(1)}\choose n_1,\ldots,n_k}
d_{\lambda^{(1)}}\cdots d_{\lambda^{(k)}}>\frac{1}{d(n^{(1)}+1)^{k(d+1)^2+k}} (a-\varepsilon)^{n^{(1)}}.
\end{equation}

There exists a multilinear polynomial
$$
f=f(x^{(1)}_{1},\ldots,x^{(1)}_{n_1},\ldots,x^{(k)}_{1},\ldots,x^{(k)}_{n_k})\not\in Id^{gr}(A)
$$
where $x^{(i)}_j\in X_{g_i}$ for all $1\le i\le k, 1\le j\le n_k$ such that $f$ generates an
irreducible $F[H]$-module with the character
$$
\chi_{_H}(F[H]f)=\chi_{_{\lambda^{(1)},\ldots,\lambda^{(k)}}}.
$$

There are homogeneous in $\Gamma$-grading
$a^{1}_{1},\ldots,a^{1}_{n_1},\ldots,a^{k}_{1},\ldots,a^{k}_{n_k}$
with $\deg_\Gamma a^i_j=g_i$ such that
$$
Q=f(a^1_1,\ldots,a^k_{n_k})\ne 0.
$$
Since $\Gamma$ is associative and commutative it follows that $Q$ is homogeneous in $\Gamma$-grading
of $A$. Therefore one can find $d^\prime\le d$ and homogeneous 
$c_1,\ldots,c_{d^\prime}\in A$ satisfying the inequality
\begin{equation}\label{eq4a}
(Q\ast c_1\ast\ldots\ast c_{d^\prime})\ast Q\ne 0
\end{equation}
where $a\ast b$ denotes the right or the left multiplication by $b$ 
(otherwise $Span<Q>$ is a graded ideal of $A$). Denote
$$
g_1=f_1=
f(x^{(1)}_{1,1},\ldots,x^{(1)}_{1,n_1},\ldots,x^{(k)}_{1,1},\ldots,x^{(k)}_{1,n_k}),
$$
where $x^{(i)}_{\alpha,\beta}$ are new homogeneous variables,
$\deg_\Gamma x^{(i)}_{\alpha,\beta}=g_i$, and take
$$
g_2=(f_1\ast z_1\ast\cdots\ast z_{d^\prime})\ast f_2,
$$
where $\deg_\Gamma z_1=\deg_\Gamma c_1,\ldots,\deg_\Gamma z_{d^\prime}=\deg_\Gamma c_{d^\prime}$,
$$
f_2=f(x^{(1)}_{2,1},\ldots,x^{(1)}_{2,n_1},\ldots,x^{(k)}_{2,1},\ldots,x^{(k)}_{2,n_k}),
$$
with $x^{(i)}_{\alpha,\beta}\in X_{g_i}$. Then $g_2\in P_{2n_1+q_1,\ldots,2n_k+q_k}$ is not
an identity of $A,$ as follows from (\ref{eq4a}), and $q_1,\ldots,q_k \ge 0, q_1+\cdots+ q_k=d^\prime$.

The square $H\times H$ of group $H$ acts on $\{x^{(i)}_{\alpha,\beta}\}$ where the first copy of 
$H$ acts on $\{x^{(i)}_{1,\beta}\}$, while the second copy acts on
$\{x^{(i)}_{2,\beta}\}$ and
$$
\chi_{_{H\times H}}(F[H\times H] g_2)=\chi_{_{\lambda^{(1)},\ldots,\lambda^{(k)}}}
\otimes \chi_{_{\lambda^{(1)},\ldots,\lambda^{(k)}}}.
$$

Denote $n^{(2)}=2n^{(1)}+d^\prime$. Repeating this procedure we construct for all $q=3,4,\ldots$ 
a multilinear polynomial
$$
g_q=g_q(x^{(1)}_{1,1},\ldots,x^{(1)}_{1,n_1},\ldots,x^{(k)}_{q,1},\ldots,x^{(k)}_{q,n_k}, z_1,z_2,\ldots)
$$
of degre $n^{(q)}$ such that:
\begin{itemize}
\item[({\it i})]
all $x^{(i)}_{\alpha,\beta}, z_j$ are homogeneous and $x^{(i)}_{\alpha,\beta}\in X_{g_i}$;
\item[({\it ii})]
$g_q$ is not an identity of $A$;
\item[({\it iii})]
$n^{(q)}=qn^{(1)}+d^{(q)},~ d^{(q)}\le (q-1)d\le dq,~ n^{(q)}-n^{(q-1)}\le d$; and
\item[({\it i}v)]
$q$ copies $H^q=H\times\cdots\times H$ of $H$ acts on $g_q$ permuting $x^{(i)}_{\alpha,\beta}$
and $g_q$ generates an irreducible $F[H^q]$-module $M$ with
$$
\chi(M)=(\chi_{_{\lambda^{(1)},\ldots,\lambda^{(k)}}})^{\otimes q}.
$$
\end{itemize}

Denote $H(q)=S_{qn_1}\times\cdots\times S_{qn_k}$. Given $1\le i\le k$, group $S_{qn_i}$
acts on $\{ x^{(i)}_{1,1},\ldots,x^{(i)}_{q,n_i}\}$. We can induce the $H^q$-action on $M$
to the $H(q)$-action. Consider the decomposition of $\widetilde M=F[H(q)] g_q$ into irreducible
components,
$$
\chi_{_H(q)}(\widetilde M)=\sum_{\rho^{(1)}\vdash qn_{1},\ldots,\rho^{(k)}\vdash qn_{k}}
t_{\rho^{(1)},\ldots,\rho^{(k)}} \chi_{_{\rho^{(1)},\ldots,\rho^{(k)}}}.
$$
It follows by the Richardson--Littlewood rule that for any $1\le i\le k$, either 
$\rho^{(i)}=q\lambda^{(i)}$ or $\rho^{(i)}$ is obtained from $\lambda^{(i)}$ by putting
down one or more boxes of $D_{q\lambda^{(i)}}$. Then by Lemma \ref{la2} we have
$\Phi(\rho^{(i)})\ge \Phi(q\lambda^{(i)})=\Phi(\lambda^{(i)})$. Now, Lemma \ref{la1} implies 
the inequality
\begin{equation}\label{eq5}
d_{\rho^{(1)}}\cdots d_{\rho^{(k)}}>\left(\frac{1}{qn^{(1)}}\right)^{k(d+1)^2} 
d_{q\lambda^{(1)}}\cdots d_{q\lambda^{(k)}}.
\end{equation}
Recall that $g_q$ is not an identity of $A$. Hence there exist integers 
$p_{q,1},\ldots, p_{q,k}\ge 0$ such that $p_{q,1}+\cdots+ p_{q,k}=d^{(k)}$ and
$$
c_{qn_1+p_{q,1},\ldots,qn_k+p_{q,k}}(A)\ge d_{\rho^{(1)}}\cdots d_{\rho^{(k)}}.
$$
In particular,
\begin{equation}\label{eq5a}
c_{n^{(q)}}^{gr}(A) \ge {qn^{(1)}\choose qn_1,\ldots,qn_k} d_{\rho^{(1)}}\cdots d_{\rho^{(k)}}.
\end{equation}

Note that for any partition $\mu\vdash m,\nu\vdash m$ with $\Phi(\mu)\ge\Phi(\nu)$ it follows by
Lemma \ref{la1} that
$$
d_\mu\ge\frac{1}{m^{d^2+d}}\Phi(\mu)^m\ge \frac{1}{m^{d^2+d}}\Phi(\nu)^m
\ge \frac{1}{m^{d^2+d+1}} d_\nu.
$$
Then by Lemma \ref{le1} and (\ref{eq5}), the right hand side of (\ref{eq5a}) is not less than
$$
\left(\frac{1}{qn^{(1)}}\right)^{2k(d+1)^2} \left(\frac{1}{n^{(1)}}\right)^{2kq}
\left[{n^{(1)}\choose n_1,\ldots,n_k} d_{\lambda^{(1)}}\cdots d_{\lambda^{(k)}}\right]^{q}.
$$
Now, (\ref{eq4}) implies the following inequality
$$
c_{n^{(q)}}^{gr}(A) \ge
\left(\frac{1}{qn^{(1)}}\right)^{2k(d+1)^2} \left(\frac{1}{n^{(1)}}\right)^{2kq}
\left(\frac{1}{d(n^{(1)}+1)}\right)^{2k(d+1)^2q} (a-\varepsilon)^{n^{(q)}}.
$$

Since $a>1$, by the assumptions of the lemma we then have 
$$
(a-\varepsilon)^{qn^{(1)}} \ge\frac{(a-\varepsilon)^{n^{(q)}}}{a^{qd}}
$$
for all small enough $\varepsilon$. Hence
$$
\sqrt[n]{c_n^{gr}(A)}>D(a-\varepsilon),
$$
where $D=D_1D_2$,
$$
D_1 =\left(\frac{1}{n^{(q)}}\right)^\frac{{2k(d+1)^2}}{n^{(q)}},\quad
D_2 =\left(\frac{1}{n^{(1)}}\right)^\frac{{2k}}{n^{(1)}}
\left(\frac{1}{d(n^{(1)}+1)}\right)^\frac{{2k(d+1)^2}}{n^{(1)}}
\left(\frac{1}{a}\right)^\frac{{d}}{n^{(1)}+d}.
$$
For small $\delta_1,\delta_2>0$ one can choose $n^{(1)}$ such that $D_1>(1-\delta_1)$ 
and $D_2>(1-\delta_2)$ for all $n^{(q)}, q\ge 1$. Finally, we can take $\delta_1,\delta_2$
small enough and get the inequality
$$
\sqrt[n]{c_n^{gr}(A)}>(1-\delta)(a-\varepsilon),
$$
for all $n=n^{(q)},~ q=1,2,\ldots~~$.
\hfill $\Box$
\begin{remark}\label{r1}
In the proof of the previous lemma we used associativity and commutativity of $\Gamma$
only for getting relation (\ref{eq4a}). In case of an arbitrary groupoid $\Gamma$ the 
element $Q$ in (\ref{eq4a}) can be non-homogeneous in $\Gamma$-grading and hence an
ideal $I$ generated by $Q$ in $A$ can be strictly less than $A$. But if $A$ is simple 
in a non-graded sense then $I=A$ and relation (\ref{eq4a}) and Lemma \ref{le2} hold.
\end{remark}

For completing the proof of the main result we need the following remark. Denote by $Ann~A$
the annihilator of $A$.
\vskip 0.2in

\begin{lemma}\label{le3}
Let $A$ be a $\Gamma$-graded algebra with a finite support of order $k$. If $Ann~A=0$ then
$$
c_{n+1}^{gr}(A)>\frac{1}{8k n^k} c_n^{gr}(A)
$$
for all sufficiently large $n$.
\end{lemma}

\pp 
Denote $Supp~A=\{g_1,\ldots,g_k\}$.  It follows from (\ref{e1}) that there exist 
$n_1,\ldots,n_k \ge 0$, $n_1+\cdots+n_k =n$ such that
\begin{equation}\label{eq6}
\frac{1}{2n^k} c_n^{gr}(A)  < \frac{1}{(n+1)^k} c_n^{gr}(A)< c_{n_1,\ldots,n_k}(A).
\end{equation}
Recall that $c_{n_1,\ldots,n_k}(A) = \dim P_{n_1,\ldots,n_k}(A)$ (see \ref{e3}). Denote
by $U_i, 1\le i\le k$, the subspace of polynomials $f$ in $P_{n_1,\ldots,n_k}(A)$ such
that $\varphi(f) A_{g_i}=0$ for all graded evaluations $\varphi: F\{X\}\to A$. Similarly,
let $W_i, 1\le i\le k$, be the subspace of polynomials $h\in P_{n_1,\ldots,n_k}(A)$ satisfying 
$A_{g_i}\varphi(h)=0$ for all graded evaluations $\varphi: F\{X\}\to A$. Denote
$$
V=U_1\cap\ldots\cap U_k\cap W_1\cap\ldots\cap W_k.
$$

If $f\in V$, then all values of $f$ in $A$ lie in $Ann~A=0$, that is $V=0$. Suppose that
$$
\dim U_1,\ldots,\dim U_k,\ldots,\dim W_1,\ldots,\dim W_k > \frac{4k-1}{4k} N
$$
where $N=c_{n_1,\ldots,n_k}(A)$. Then $\dim V > (2k\cdot\frac{4k-1}{4k}-(2k-1))N=\frac{1}{2}N$,
that is $V\ne 0$, a contradiction. It follows that $\dim U_i < \frac{4k-1}{4k} N$ or
$\dim W_i < \frac{4k-1}{4k} N$ for at least one $i$. Let, for instance, $\dim U_1 < \frac{4k-1}{4k} N$.
Denote by $T$ the codimension of $U_1$ in $P_{n_1,\ldots,n_k}(A)$. Then
$$
T > \frac{1}{4k} c_{n_1,\ldots,n_k}(A) > \frac{1}{8k n^k} c_n^{gr}(A)
$$
as follows from (\ref{eq6}). Now if $f_1,\ldots, f_T$ are linearly independent modulo $U_1$ elements 
from $P_{n_1,\ldots,n_k}(A)$ then $f_1z,\ldots, f_Tz$ are linearly independent elements in
$P_{n_1+1,n_2,\ldots,n_k}(A)$, provided that $z$ is a new homogeneous indeterminate, $\deg_\Gamma z= g_1$.
Hence
$$
c_{n+1}^{gr}(A) \ge  c_{n_1+1,n_2,\ldots,n_k}(A) \ge
T >  \frac{1}{8k n^k} c_n^{gr}(A),
$$
and we are done.
\hfill $\Box$
\vskip 0.2in

Now we are ready to prove the main result of this paper.

\begin{theorem}\label{t2}
Let $\Gamma$ be a commutative semigroup and let $A=\bigoplus_{g\in\Gamma} A_g$ be a finite dimensional 
$\Gamma$-graded algebra. If $A$ is graded simple then there exists the limit
$$
exp^{gr}(A)=\lim_{n\to\infty} \sqrt[n]{c_n^{gr}(A)}.
$$
\end{theorem}

\pp
Denote
$$
a=\overline{exp}^{gr}(A)=\limsup_{n\to\infty} \sqrt[n]{c_n^{gr}(A)}.
$$
If $a=0$ then $A$ is nilpotent and $exp^{gr}(A)=0$. If $A$ is not nilpotent then $a\ge 1$. In the case $a=1$ the
lower limit of $\sqrt[n]{c_n^{gr}(A)}$ is also 1 and we are done.

Let now $a>1$. By Lemma \ref{le2} there exists a sequence $n^{(1)}< n^{(2)}<\ldots$ such  that
$$
c_n^{gr} \ge (1-\delta)^n(a-\varepsilon)^n
$$
for all $n=n^{(i)}, i\ge 1$, and $\varepsilon,\delta>0$ can be choosen arbitrary small.

Now let $m=n^{(i)}, m^\prime =n^{(i+1)}$ and let $m<n<m^\prime$. Then $n=m+p, 1\le p<d$. By Lemma \ref{le3}
we have
\begin{equation}\label{eq7}
c_n^{gr}(A)= c_{m+p}^{gr}(A) > \left( \frac{1}{(8k(m+p))}\right)^p (1-\delta)^n(a-\varepsilon)^m
>\left(\frac{1}{8kn}\right)^d \frac{1}{(a-\varepsilon)^d} (1-\delta)^n(a-\varepsilon)^n.
\end{equation}

Clearly, inequality (\ref{eq7}) also holds  for all $n=n^{(1)}, n^{(2)},\ldots$, and for all small
$\varepsilon,\delta>0$. Hence
$$
\liminf_{n\to\infty} \sqrt[n]{c_n^{gr}(A)} \ge (1-\delta)(a-\varepsilon).
$$
Since $\varepsilon,\delta$ are arbitrary we can conclude that
$$
\underline{exp}^{gr}(A)=\liminf_{n\to\infty} \sqrt[n]{c_n^{gr}(A)}=a
$$
and the proof of the theorem is completed.
\hfill $\Box$

Finally, note that associativity and commutativity of $\Gamma$ was used only in the proof of Lemma \ref{le2}
(see Remark \ref{r1}). Hence for an arbitrary groupoid $\Gamma$ we have obtained the following result.
\vskip 0.2in

\begin{theorem}\label{t3}
Let  $A=\bigoplus_{g\in\Gamma} A_g$ be a finite dimensional algebra graded by a groupoid $\Gamma$. 
If $A$ is simple then there exists its graded PI-exponent
$$
exp^{gr}(A)=\lim_{n\to\infty} \sqrt[n]{c_n^{gr}(A)}.
$$
\end{theorem}

\hfill $\Box$

\end{document}